\documentclass[11pt]{amsart}

\usepackage{amssymb}

\usepackage{enumerate}

\usepackage{graphicx}

\makeatletter
\@namedef{subjclassname@2020}{%
  \textup{2020} Mathematics Subject Classification}
\makeatother

\newtheorem{thm}{Theorem}[section]
\newtheorem{cor}[thm]{Corollary}
\newtheorem{lem}[thm]{Lemma}
\newtheorem{prop}[thm]{Proposition}

\theoremstyle{definition}
\newtheorem{defn}[thm]{Definition}
\newtheorem{rem}[thm]{Remark}
\newtheorem{exa}[thm]{Example}

\numberwithin{equation}{section}

\begin{document}

\baselineskip=17pt

\title[Dynamics of unbounded linear operators]{Dynamics of unbounded linear operators}

\author[Mohammad Ansari]{Mohammad Ansari}

\date{}
\begin{abstract} We apply the well-known and also the newly introduced notions from bounded linear dynamics to unbounded linear operators.
We present a hypermixing criterion similar to that given for bounded linear operators and
then we show that the derivative operator in $H^2$, the Laplacian operator in $L^2(\Omega)$, 
 and all unbounded weighted translations in $L_p(0,\infty)$ and $C_0[0,\infty)$ are hypermixing.
\end{abstract}

\maketitle
\section{Introduction}
Motivated by the papers \cite{bcs} and \cite{jm} in which the authors investigate the chaoticity of some unbounded linear operators,
 we decided to apply the other well-known and also the newly introduced notions from bounded linear dynamics to unbounded linear operators. 
\par Perhaps the most well-known dynamical notion which has been extensively investigated for bounded linear operators is {\it hypercyclicity}.
 Let $\mathcal X$ be a separable infinite-dimensional (real or complex) Banach space and let $T: \mathcal X\to \mathcal X$
be a bounded linear operator. Then $T$ is said to be {\it hypercyclic} if there exists a vector $x\in\mathcal X$ such that the orbit of $x$ under $T$
$$\textnormal{orb}(x, T) = \{x, Tx, T^2x, T^3x,\cdots\}$$ is dense in $\mathcal X$, and any such vector $x$ is called a {\it hypercyclic
vector} of $T$.
In \cite{r}, Rolewicz offered the first example of a hypercyclic bounded linear operator on a Banach space: The operator $\lambda B$ where $B$ is the unilateral backward
shift on $\ell^p$ $(1\le p<\infty)$ or $c_0$, and $|\lambda|>1$. 
\par It is clear that if $x$ is a hypercyclic vector for $T$, then so is every $T^px$ ($p\ge 1$), and hence, for the set $HC(T)$ of all hypercyclic vectors for $T$, 
  we have that $HC(T)=\emptyset$ or $\overline{HC(T)}=\mathcal X$.
\par Let us also recall the notion of {\it chaoticity} for bounded linear operators. If $T^Nx=x$ for some positive integer $N$ and some $x\in \mathcal X$, 
then $x$ is called a {\it periodic vector} for $T$. The subspace of $\mathcal X$ comprising of all periodic vectors for $T$ is denoted by $\textnormal{Per} (T)$. The operator
$T$ is said to be {\it chaotic} (following the definition of
Devaney \cite{d}) if it is hypercyclic and has a dense set of periodic vectors, or equivalently, if $\overline{HC(T)}=\mathcal X=\overline{\textnormal{Per} (T)}$.
For more information on hypercyclic and chaotic operators, the reader is referred to \cite{bm} and \cite{gp}.
\par The notions of hypercyclicity and chaoticity can be similarly defined for unbounded linear operators, but one should consider domain-related conditions.
More precisely, if $T$ is an unbounded linear operator in $\mathcal X$ with domain $D(T)$, then it is clear that if we want to talk about the orbit of a vector $x\in \mathcal X$ 
under $T$, we need to have $x\in \bigcap_{n=1}^{\infty}D(T^n)$. Thus $T$ is hypercyclic if there is some $x\in \mathcal X$ such that 
$$x\in \bigcap_{n=1}^{\infty}D(T^n)\;\textnormal{and}\; \overline{\textnormal{orb}(x, T)}=\mathcal X.$$
 On the other hand, to be a periodic vector for $T$, a vector $x\in \mathcal X$ should satisfy $$x\in D(T^N)\; \textnormal{and}\; T^Nx=x,$$
for some positive integer $N$. Here too, we have $HC(T)=\emptyset$ or $\overline{HC(T)}=\mathcal X$ (because $x\in \bigcap_{n=1}^{\infty}D(T^n)$
implies $T^px\in \bigcap_{n=1}^{\infty}D(T^n)$ for all $p\ge 1$), but it may not seem so obvious that $\textnormal{Per} (T)$ is a linear subspace
of $\mathcal X$. In fact, it may not even be seen at a glance that $\textnormal{Per} (T)\subseteq Y_0$. 
\begin{lem} Let $T$ be an unbounded linear operator in $\mathcal X$. Then $\textnormal{Per} (T)\subseteq Y_0$, and as a result, $\textnormal{Per} (T)$ 
is a linear subspace of $\mathcal X$.
\begin{proof} Let $x\in\textnormal{Per} (T)$ be arbitrary. Then there is some positive integer $N$ such that $x\in D(T^N)$ and $T^Nx=x$.
Then it is clear that $x\in D(T^{mN})$ for all integers $m\ge 1$. Now assume that $n\ge 1$ is an arbitrary integer. Then there is some $m\ge 1$ such that $mN\ge n$, and hence,
$x\in D(T^{mN})\subseteq D(T^n)$. Thus we have that $x\in Y_0$ which proves that $\textnormal{Per} (T)\subseteq Y_0$ since $x$ was arbitrary.
Now we prove that $\textnormal{Per} (T)$ is a linear subspace. It is clear that $x\in \textnormal{Per} (T)$ implies that $\lambda x\in \textnormal{Per} (T)$ for any scalar $\lambda$.
On the other hand, if $x, y\in\textnormal{Per} (T)$ then $x, y\in Y_0$, and hence, $x+y\in Y_0$. Now, if $T^Nx=x$ and $T^Ly=y$ for some positive integers $N, L$, then $T^{NL}(x+y)=x+y$.
Thus $x+y\in\textnormal{Per} (T)$.
\end{proof}
\end{lem}
 In \cite{bcs} the authors provide the following sufficient condition (analogous to the hypercyclicity
criterion for continuous linear operators obtained by Kitai \cite{k} and independently by Gethner and Shapiro \cite{gs}) for a (bounded or unbounded) linear 
operator to be hypercyclic. Recall that a linear operator $T$ is said to be {\it densely defined} in $\mathcal X$ if $D(T)$ is dense in $\mathcal X$.
 Moreover, $T$ is called {\it closed} if 
its graph $\{(x,Tx):x\in D(T)\}$ is closed in $X\times X$.
\begin{thm} [Theorem $2.1$ of \cite{bcs}]  Suppose that $\mathcal X$ is a separable infinite-dimensional 
 Banach space and let $T$ be a densely defined linear operator in $\mathcal X$ for which any power $T^n$ ($n\ge 1$)
is a closed operator. Then $T$ is hypercyclic if there exists a set $Y\subseteq \bigcap_{n=1}^{\infty}D(T^n)$ dense in $X$ and a mapping $S:Y\to Y$ such that \\
 \hspace*{.8cm} \textnormal{i)} $TSx=x$ for all $x\in Y$, and\\
\hspace*{.65cm} \textnormal{ii)} $T^nx$, $S^nx\to 0$ for all $x\in Y$.
\end{thm}
Then they apply the condition to the
derivative operator in $H^2$ and the Laplacian operator in $L^2(\Omega)$ to show that these operators are hypercyclic. Also, by proving the existence of dense
sets of periodic vectors, the authors show that the mentioned unbounded operators are chaotic. 
\par On the other hand, in their recently published paper \cite{jm}, the authors use Theorem $1.2$ above to 
show that all unbounded weighted translation operators in $L_p(0,\infty)$ and $C_0[0,\infty)$ are 
hypercyclic, and then they prove that the mentioned operators are chaotic by showing that they possess dense sets of periodic vectors.
\par In Section $2$, motivated by the domain restriction of unbounded linear operators, 
we present the modified versions of the notions of topological transitivity (TT), weakly mixing (WM), mixing (M), strong topological transitivity (STT),
 supermixing (SM), and hypermixing (HM) for linear operators. Also, we offer the notion of {\it relative hypercyclicity} (RHC) and we give some results concerning its relationships with 
 hypercycliciity (HC) and topological transitivity. Then we 
state two criteria for strong topological transitivity and hypermixing, analogous to those given for bounded linear operators \cite[Theorems $2.1$ and $2.3$]{a2}.
Finally, we offer the notion of {\it up-to-closure extension} for a dynamical domain and we give some relevant results. 
\par In Section $3$, we prove that the derivative operator in the Hardy Hilbert space $H^2$, the Laplacian operator in $L^2(\Omega)$, 
 and all unbounded weighted left translations in $L_p(0,\infty)$ ($1\le p<\infty$) and $C_0[0,\infty)$ are hypermixing.
\section{Dynamical notions with respect to a subspace}
As we have mentioned in the introduction, the notions of hypercyclicity and chaoticity have been studied for unbounded linear operators. In this 
section, we present other dynamical notions for these operators in a form which is also applicable to bounded linear operators.
\par Recall that a bounded linear operator $T: \mathcal X\to \mathcal X$ is called {\it topologically transitive} if, for any pair of nonempty
 open subsets $U, V$ of $\mathcal X$, there is a non-negative integer $n$ such that $T^n(U)\cap V\neq \emptyset$. 
 If $T\oplus T$ is topologically transitive on $\mathcal X\oplus \mathcal X$, we say that $T$ is {\it weakly mixing}.
Finally, $T$ is said to be {\it mixing} whenever, for any pair of nonempty
 open subsets $U, V$ of $\mathcal X$, there exists a non-negative integer $N$ such that $T^n(U)\cap V\neq\emptyset$ for all $n\ge N$.
 \par Now we consider the above-mentioned dynamical concepts for unbounded linear operators.
 Let $T: D(T)\subseteq \mathcal X\to \mathcal X$ be an unbounded linear operator. Since the concepts involve the effects 
 of the powers $T^n$ on open sets, we clearly must work with those nonempty sets $U$ which are (relatively) open in the subspace $Y_0=\bigcap_{n=1}^{\infty}D(T^n)$.
 In fact, we offer our definitions in a more general setting: For
 a linear (not necessarily closed) subspace $Y$ of $\mathcal X$ satisfying $Y\subseteq Y_0$, we define the dynamical concepts for $T$ by considering the nonempty (relatively) open sets in $Y$,
 and we say that $Y$ is the underlying {\it dynamical domain} of $T$. 
 Then $Y=Y_0$ would be a particular case and we call $Y_0$ the {\it maximal dynamical domain} of $T$. We use the notation $(T,Y)$ to refer to the subspace $Y$, and
clearly it makes sense if we write $T$ for $(T,Y_0)$ and say that $T$ has a dynamical property whenever $(T,Y_0)$ has that property. 
Note that the notion of hypercyclicity is irrelevant to the underlying dynamical domain, and hence, we never say that $(T, Y)$ is hypercyclic.
 \par The usefulness of considering the general case $(T,Y)$ is that we can present our definitions and results for both unbounded and bounded linear operators.
 If we only work with the maximal dynamical domain $Y_0$, then it gives the usual dynamical notions for bounded linear operators (because $Y_0=\mathcal X$ for bounded operators) 
 which have already been investigated.
\par Therefore, throughout this paper,\\
 $\bullet$ $\mathcal X$ is a separable infinite-dimensional (real or complex) Banach space,\\
 $\bullet$ $T: D(T)\subseteq \mathcal X\to \mathcal X$ is a linear operator,\\
 $\bullet$ $Y_0=\bigcap_{n=1}^{\infty}D(T^n)$ is the maximal dynamical domain of $T$ and
 $Y$ is a linear subspace of $\mathcal X$ such that $Y\subseteq Y_0$ (or equivalently, $Y$ is a linear subspace of $Y_0$),\\
 $\bullet$ We write $(T,Y)$ to say that $Y$ is the underlying dynamical domain of $T$,\\
 $\bullet$ If $(T,Y_0)$ has a dynamical property, we say that $T$ has that property and we write $T$ instead of $(T,Y_0)$, and\\
 $\bullet$ For any set $M\subseteq Y$, by $\overline M$ we mean the closure of $M$ in $\mathcal X$ (we can also work with the closures in $Y$, but since we will have connections to 
$\overline{\textnormal{orb}(x, T)}$ which is a closure in $\mathcal X$, we choose to consider the closures in $\mathcal X$).
\par Now we are ready to apply the notions of topological transitivity, weakly mixing, and mixing to the pair $(T,Y)$.
 \begin{defn} Suppose that $T$ is a linear operator in $\mathcal X$. We say that $(T,Y)$ is\\
\hspace*{.2cm} \textnormal{i)} {\it topologically transitive} if, for any pair of nonempty open sets $U, V$ in $Y$, we have $T^n(U)\cap V\neq \emptyset$ for
some non-negative integer $n$,\\
\hspace*{.1cm} \textnormal{ii)} {\it weakly mixing} whenever $(T\oplus T, Y\oplus Y)$ is a topologically transitive operator in $\mathcal X\oplus \mathcal X$, and\\
\hspace*{.03cm} \textnormal{iii)} {\it mixing} if, for any pair of nonempty
 open sets $U, V$ in $Y$, there exists a non-negative integer $N$ such that $T^n(U)\cap V\neq \emptyset$ for all $n\ge N$.
 \end{defn}
 Recently, some stronger versions of topological transitivity and mixing have been introduced and investigated \cite{a1,a2,a3,a4,c}.
 A bounded linear operator $T$ on $\mathcal X$ is said to be {\it strongly topologically transitive} if, for any nonempty open subset $U$ of $\mathcal X$, we have
 $\mathcal X\backslash\{0\}\subseteq \bigcup_{n=0}^{\infty}T^n(U)$. We say that $T$ is {\it supermixing} if 
 $\mathcal X=\overline{\bigcup_{i=0}^{\infty}\bigcap_{n=i}^{\infty}T^n(U)}$ for each nonempty open set $U$ in $\mathcal X$.
 Finally, $T$ is called {\it hypermixing} whenever, for any nonempty open subset $U$ of $\mathcal X$, we have that $\mathcal X\backslash\{0\}\subseteq \bigcup_{i=0}^{\infty}\bigcap_{n=i}^{\infty}T^n(U)$, or equivalently, $\mathcal X=\bigcup_{i=0}^{\infty}\bigcap_{n=i}^{\infty}T^n(U)$ \cite[Remark $2.5$]{a2}.
 \begin{defn} Suppose that $T$ is a linear operator in $\mathcal X$. We say that $(T,Y)$ is\\
\hspace*{.2cm} \textnormal{i)} {\it strongly topologically transitive} if $Y\backslash\{0\}\subseteq \bigcup_{n=0}^{\infty}T^n(U)$,\\
\hspace*{.1cm} \textnormal{ii)} {\it supermixing} if $\overline Y=\overline{\big{(}\bigcup_{i=0}^{\infty}\bigcap_{n=i}^{\infty}T^n(U)\big{)}\cap Y}$, and\\
\hspace*{.02cm} \textnormal{iii)} {\it hypermixing} provided that $Y\backslash\{0\}\subseteq \bigcup_{i=0}^{\infty}\bigcap_{n=i}^{\infty}T^n(U)$,
\\for any nonempty open set $U$ in $Y$.
 \end{defn}
 \begin{rem} As the reader has noticed, $(T,Y)$ is supermixing if and only if, $\big{(}\bigcup_{i=0}^{\infty}\bigcap_{n=i}^{\infty}T^n(U)\big{)}\cap Y$ is dense in $Y$.
 Also, we can easily show that if $(T,Y)$ is hypermixing then $T$ is not injective, and as a result,
the definition of hypermixing for $(T,Y)$ can be modified to $Y\subseteq \bigcup_{i=0}^{\infty}\bigcap_{n=i}^{\infty}T^n(U)$. To prove these assertions, 
pick $0\neq y\in Y$
and let $U\neq \emptyset$ be open in $Y$. Let $V$ be another nonempty open set in $Y$ such that
 $U\cap V=\emptyset$. Then, by the definition of hypermixing, we can find a positive integer $N$ and two vectors $u\in U$ and $v\in V$ such that 
 $T^Nu=y=T^Nv$. This shows that $T$ is not injective. Now, if we put $x=u-v$ then $0\neq x\in Y$, and hence, there is some $N_0\ge 0$ such that $x\in T^n(U)$ for all $n\ge N_0$.
 Then $0=T^Nx\in T^n(U)$ for all $n\ge N+N_0$ which shows that $0\in \bigcup_{i=0}^{\infty}\bigcap_{n=i}^{\infty}T^n(U)$.
 Since $U\neq \emptyset$ was arbitrarily chosen, this proves that we can modify the definition of hypermixing to $Y\subseteq \bigcup_{i=0}^{\infty}\bigcap_{n=i}^{\infty}T^n(U)$ for all
 $U\neq \emptyset$ open in $Y$.
 \end{rem}
 The following diagram of implications for the pair $(T,Y)$ follows quickly from Definitions $2.1$ and $2.2$:
 \\
\[
\begin{array}
[c]{cccccccc}
   \text{HM}&& \Rightarrow& &\text{SM}\\
         \vspace*{.2cm} && && \Downarrow&\\
    \Downarrow && &&\text{M}\\
   && && \Downarrow\\
   \text{STT} &\Rightarrow&\text{TT}&\Leftarrow& \text{WM}\vspace*{.2cm}
\end{array}
\]
\\
\par Recall that, by the {\it Birkhoff transitivity theorem} \cite{bi}, topological transitivity and hypercyclicity
are equivalent ($TT\Leftrightarrow HC$) for a bounded linear operator. 
But we may not have the same equivalence for the general case $(T,Y)$, because the density of $\textnormal{orb} (x, T)$ in $\mathcal X$
implies that the orbit intersects all nonempty open sets in $\mathcal X$, and there is no guarantee that the orbit also intersects all nonempty open 
sets in $Y$.
\par We see that if $Y$ includes the orbit of a hypercyclic vector for $T$ then $(T,Y)$ is topologically transitive. The particular cases are 
$$Y=Y_x=\textnormal{span} \big{(}\textnormal{orb}(x,T)\big{)},$$ and
$Y=Y_0$.
 \begin{prop} Assume that $T$ is a hypercyclic linear operator in $\mathcal X$. If, for some $x\in HC(T)$, we have $\textnormal{orb} (x, T)\subseteq Y$
 then $(T,Y)$ is topologically transitive. In particular, $(T, Y_x)$ is topologically transitive for all $x\in HC(T)$, and also $T$ is topologically transitive.
 \begin{proof} Suppose $x\in HC(T)$ and $\textnormal{orb} (x, T)\subseteq Y$. Let $U, V$ be nonempty open sets in $Y$ and $G, H$ be open sets in $\mathcal X$ 
 such that $U=G\cap Y$ and $V=H\cap Y$.
 Then there is a non-negative integer $p$ such that $$T^px\in G\cap \textnormal{orb}(x,T)\subseteq G\cap Y=U.$$
 Also, since $T^px\in HC(T)$, there exists a non-negative integer $q$ such that $$T^q(T^px)\in H\cap \textnormal{orb}(T^px,T)\subseteq H\cap Y=V,$$
 and hence, $T^q(U)\cap V\neq\emptyset$.
 \par To prove the particular cases, it is clear that $\textnormal{orb} (x, T)\subseteq Y_x$ and $\textnormal{orb} (x, T)\subseteq Y_0$
 for every $x\in HC(T)$.
 Thus, by the proof of the first part, $(T, Y_x)$ is topologically transitive for every $x\in HC(T)$, and also $T=(T, Y_0)$ is topologically transitive.
 \end{proof}
 \end{prop}
  Thus, in view of the particular case $Y=Y_0$ in the above proposition, the implication $HC\Rightarrow TT$ is always true for a linear operator $T$ in $\mathcal X$.
 \par Now we define the notion of 
{\it relative hypercyclicity} for $(T,Y)$ in order to make a connection to topological transitivity of $(T,Y)$. Before it, let us recall the notion of {\it subspace hypercyclicity}
for bounded linear operators \cite{mm}: Let $\mathcal M$ be a closed linear subspace of $\mathcal X$. A bounded linear operator $T$ on $\mathcal X$ is said to be {\it subspace hypercyclic}
for $\mathcal M$ if $\textnormal{orb}(x,T)\cap \mathcal M$ is dense in $\mathcal M$, or equivalently 
(using the closure in $\mathcal X$), $\overline{\textnormal{orb}(x,T)\cap \mathcal M}\supseteq \mathcal M$. Clearly, this notion can be similarly defined for an unbounded operator.
 \par The authors in \cite{re} consider the notion of subspace hypercyclicity for general subspaces (which may not necessarily be closed)
 based on their claim that ``a bounded linear operator $T$ on $\mathcal X$ is subspace hypercyclic for $\mathcal M$ if and only if $T$ is subspace hypercyclic for $\overline{\mathcal M}$''.
 But this claim is not true (only the necessity condition is true), because if $T$ is a chaotic bounded linear operator on $\mathcal X$ and $\mathcal M=\textnormal{Per}(T)$, then we can observe that $T$ is subspace
 hypercyclic for $\overline{\mathcal M}=\mathcal X$ but $T$ is not subspace hypercyclic for $\mathcal M$ (see Example $2.10$ (ii)).
\par Therefore, we choose a different name for the following similar notion 
which is presented for general linear subspaces.
\begin{defn} Suppose $T$ is a linear operator in $\mathcal X$. Then the pair $(T,Y)$ is called
{\it relatively hypercyclic} if there exists a nonzero vector $x\in Y_0$ for which $\overline{\textnormal{orb}(x,T)\cap Y}\supseteq Y$.
Such a vector $x$ is called a {\it relatively hypercyclic vector} for $(T,Y)$, and the set of all relatively hypercyclic vectors for $(T,Y)$ is denoted by $RHC(T,Y)$.
We say that $T$ is relatively hypercyclic whenever $(T,Y_0)$ is relatively hypercyclic and, in that case, we write $RHC(T)$ for $RHC(T,Y)$.
\end{defn}
\begin{rem} It is clear that the two notions of relative hypercyclicity and subspace hypercyclicity 
coincide whenever $Y$ is a closed subspace of $\mathcal X$. Meanwhile, note that 
the condition $\overline{\textnormal{orb}(x,T)\cap Y}\supseteq Y$ is equivalent to $\overline{\textnormal{orb}(x,T)\cap Y}=\overline{Y}$.
Also, we have observed that $HC(T)=\emptyset$ or $\overline{HC(T)}=\mathcal X$. Now, for the set 
$RHC(T,Y)$, we claim that $RHC(T,Y)=\emptyset$ or $\overline{RHC(T,Y)\cap Y}=\overline Y$ (or equivalently, $RHC(T,Y)\cap Y$ is dense in $Y$). Indeed, if $x\in RHC(T,Y)$
then so does every $T^px$ ($p\ge 1$) because $$\textnormal{orb}(T^px,T)\cap Y=(\textnormal{orb}(x,T)\cap Y)\backslash{\{T^nx: 0\le n\le p-1\}}$$ 
and $\overline{\textnormal{orb}(x,T)\cap Y}=\overline{Y}$.
This shows that $\textnormal{orb}(x,T)\subseteq RHC(T,Y)$. Then $$\overline Y=\overline{\textnormal{orb}(x,T)\cap Y}\subseteq \overline{RHC(T,Y)\cap Y}\subseteq \overline Y.$$
\end{rem}
In the next three results, we deal with the possibility of implications between relative hypercyclicity on one hand, and topological transitivity and hypercyclicity on the other.
\begin{prop} Suppose that $T$ is a linear operator in $\mathcal X$.
If $(T,Y)$ is relatively hypercyclic then $(T,Y)$ is topologically transitive.
\begin{proof} Let $U, V\neq\emptyset$ be open sets in $Y$. Since $(T,Y)$ is relatively hypercyclic, there is some $x\in Y_0$
so that $\overline{\textnormal{orb}(x,T)\cap Y}\supseteq Y$. Let $G$ be an open set in $\mathcal X$ for which $U=G\cap Y$.
If we pick some $u\in U$ then $u\in \overline{\textnormal{orb}(x,T)\cap Y}$. But $u\in G$ and $G$ is open in $\mathcal X$, and hence, $G$ should intersect
$\textnormal{orb}(x,T)\cap Y$. Then there exists a non-negative integer $p$ such that $T^px\in G\cap Y=U$. On the other hand, by remark $2.6$, $T^px$ is also a relatively hypercyclic 
vector for $(T,Y)$, and hence, $\overline{\textnormal{orb}(T^px,T)\cap Y}\supseteq Y$. Now, for the open set $V$, there is an open set $H$ in $\mathcal X$
such that $V=H\cap Y$, and similar to what we have mentioned for the set $G$, the set $H$ should  intersect $\textnormal{orb}(T^px,T)\cap Y$. 
Then there is a non-negative integer $q$ such that $T^q(T^px)\in H\cap Y=V$.
Hence, $T^q(U)\cap V\neq \emptyset$ and we are done.
\end{proof}
\end{prop}
Hence we always have $RHC\Rightarrow TT$ for the pair $(T,Y)$.
\begin{prop} Let $T$ be a linear operator in $\mathcal X$. If $\overline Y=\mathcal X$
and $(T,Y)$ is relatively hypercyclic then $T$ is hypercyclic, and moreover, $RHC(T,Y)\subseteq HC(T)$.
\begin{proof} Since $(T,Y)$ is relatively hypercyclic, we have that $RHC(T,Y)\neq\emptyset$. Pick an arbitrary $x\in RHC(T,Y)$. Then
 $$\overline{\textnormal{orb}(x,T)}\supseteq\overline{\textnormal{orb}(x,T)\cap Y}=\overline{Y}=\mathcal X,$$ and hence,
  we conclude that $T$ is hypercyclic and $RHC(T,Y)\subseteq HC(T)$.
\end{proof}
\end{prop}
 Thus, if $\overline Y=\mathcal X$ then, by Propositions $2.7$ and $2.8$, for the pair $(T,Y)$ we have that
  $$HC\Leftarrow RHC\Rightarrow TT.$$
  \par Recall that, if a bounded linear operator is hypercyclic then it is subspace hypercyclic \cite{bkk}. The following result gives a similar implication between hypercyclicity and relative 
  hypercyclicity of an unbounded linear operator. In fact, it also gives something more.
  Since for a bounded linear operator $T$ on $\mathcal X$ we have that $Y_0=\mathcal X$, and also the relative hypercyclicity of $(T,\mathcal X)$ 
   coincides with the hypercyclicity of $T$, we observe that the statement of the following proposition is useless for bounded linear operators.
\begin{prop} Suppose that $T$ is an unbounded linear operator in $\mathcal X$. Then $T$ is hypercyclic if and only if $T$ is relatively hypercyclic
and $\overline {Y_0}=\mathcal X$. In that case, $RHC(T)=HC(T)$.
\begin{proof} The sufficiency has already been proved in Proposition $2.8$. To prove the necessity, 
suppose $T$ is hypercyclic and let $x\in HC(T)$ be arbitrary. Then
$$\overline {Y_0}\supseteq \overline{\textnormal{orb}(x,T)\cap Y_0}=\overline{\textnormal{orb}(x,T)}=\mathcal X\supseteq Y_0,$$
which shows that $T$ is relatively hypercyclic, $HC(T)\subseteq RHC(T)$, and $\overline {Y_0}=\mathcal X$. Now, by Proposition $2.8$, we also have 
$RHC(T)\subseteq HC(T)$, and hence, $RHC(T)=HC(T)$. 
\end{proof}
\end{prop}
Hence, by Propositions $2.7$ and $2.9$, we have the following chain of implications for $T$:
$$HC\Rightarrow RHC\Rightarrow TT.$$
If moreover, $\overline {Y_0}=\mathcal X$ then, by Propositions $2.7$ and $2.9$, we have that
 $$HC\Leftrightarrow RHC\Rightarrow TT.$$
 \par In the following example, first we give a relatively hypercyclic operator which is not hypercyclic, and then
 we present a hypercyclic operator $T$ and some subspace $Y\subseteq Y_0$ such that $(T,Y)$ is not relatively hypercyclic.
 \begin{exa} (i) Let $A=\{a_1, a_2, a_3, \cdots\}$ be a countable dense subset of $\mathcal X$ and $Y=\textnormal{span} (A\oplus (0))$, where by $(0)$ 
we mean the null subspace of $\mathcal X$. Define the operator $T$ in $\mathcal X\oplus \mathcal X$ with $D(T)=Y$ linearly 
by $T(a_i,0)=(a_{i+1},0)$ ($i\ge 1$). Then $Y_0=Y$ and
 $$\overline{\textnormal{orb}((a_1,0),T)\cap Y_0}=\overline{\textnormal{orb}((a_1,0),T)}=\mathcal X\oplus (0)\supseteq Y_0$$
which proves that $T$ is relatively hypercyclic. On the other hand, $\overline {Y_0}\neq \mathcal X\oplus \mathcal X$, and hence,
we conclude that $T$ is not hypercyclic by Proposition $2.9$.
\par (ii) Suppose $T$ is a chaotic linear operator in $\mathcal X$. If we put $Y=\textnormal{Per} (T)$ 
then ($Y\subseteq Y_0$ is a linear subspace by Lemma $1.1$ and) we claim that, for any $x\in Y_0$, 
$\overline{\textnormal{orb}(x,T)\cap Y}\neq \overline Y$.  
Indeed, if $x\in HC(T)$ then $\overline{\textnormal{orb}(x,T)\cap Y}=\emptyset$, and whenever $x\notin HC(T)$
then $\overline{\textnormal{orb}(x,T)\cap Y}\subseteq \overline{\textnormal{orb}(x,T)}\subsetneq \mathcal X=\overline Y$. Thus $(T,Y)$ is not relatively hypercyclic.
\end{exa}
The next result gives information on the impact of dynamical properties on the range of a linear operator. Recall that by $R(T)$ we mean the range of $T$.
 \begin{prop} Let $T$ be a linear operator in $\mathcal X$. Then the
 following hold.\\
\hspace*{.2cm} \textnormal{i)} If $(T,Y)$ is strongly topologically transitive then $Y\subseteq R(T)$.\\
\hspace*{.08cm} \textnormal{ii)} If $(T,Y)$ is supermixing or relatively hypercyclic then $\overline Y\subseteq\overline{R(T)}$.\\
\hspace*{.05cm}\textnormal{iii)} If $\overline Y=\mathcal X$ and $(T,Y)$ has either of the properties in \textnormal{(i)}, \textnormal{(ii)}
 then $T$ has a dense range, i.e. $\overline{R(T)}=\mathcal X$.\\
\hspace*{.05cm}\textnormal{iv)} If $T$ is hypercyclic then $\overline{R(T)}=\mathcal X$.
\begin{proof}  \textnormal{i)} Suppose $U$ is a nonempty open subset of $Y$ and let $0\neq y\in Y$ be arbitrary. Then, by the definition of strong topological transitivity for $(T,Y)$, we
 have $y\in\bigcup_{n=0}^{\infty}T^n(U)$.
 If $y\notin U$ then $y\in T^N(U)$ for some $N\ge 1$, and hence, $y\in R(T)$. Now
assume that $y\in U$. Then $V=U\backslash\{y\}$ is a nonempty open set in $Y$ and $y\notin V$, and since $Y\subseteq \bigcup_{n=0}^{\infty}T^n(V)$ we must have
$y\in T^k(V)$ for some $k\ge 1$, which says that $y\in R(T)$.\\
\hspace*{.5cm}  \textnormal{ii)} Suppose $(T,Y)$ is supermixing and let $U$ be a nonempty open set in $Y$. Then
$\overline Y=\overline{(\bigcup_{i=0}^{\infty}\bigcap_{n=i}^{\infty}T^n(U))\cap Y}\subseteq \overline{R(T)\cap Y}\subseteq \overline{R(T)}$. 
\par Now assume that $(T,Y)$ is relatively hypercyclic
 and let $x\in RHC(T,Y)$. Then, in view of Remark $2.6$, $Tx\in RHC(T,Y)$ and hence,
 $\overline Y=\overline{\textnormal{orb}(Tx,T)\cap Y}\subseteq \overline{\textnormal{orb}(Tx,T)}\subseteq \overline{R(T)}.$ \\
\hspace*{.4cm}  \textnormal{iii)} This is trivial.\\
\hspace*{.4cm}  \textnormal{iv)} Suppose that $T$ is hypercyclic and let $x\in HC(T)$. Then $Tx\in HC(T)$, and hence, $\mathcal X=\overline{\textnormal{orb}(Tx,T)}\subseteq \overline{R(T)}$.
\end{proof}
\end{prop}
By the {\it generalized kernel} of a linear operator $T$ in $\mathcal X$, we mean $$\textnormal{GK}(T)=\bigcup_{n=1}^{\infty}\text{Ker}T^n,$$
 where $\text{Ker}T^n=\{x\in \mathcal X: T^nx=0\}$ ($n\ge 1$).
\begin{prop} Suppose $T$ is a linear operator in $\mathcal X$. If $(T,Y)$ is strongly topologically transitive
then either $\textnormal{GK}(T)\cap Y=(0)$ or $\overline{\textnormal{GK}(T)\cap Y}\supseteq Y$.
\begin{proof} If $T$ is injective, then $\textnormal{GK}(T)\cap Y=(0)$. Suppose $T$ is not injective and let $y\in Y$ be arbitrary.
If $y\in \textnormal{GK}(T)$ then $y\in\overline{\textnormal{GK}(T)\cap Y}$. Assume that $y\notin \textnormal{GK}(T)$.
Pick a countable neighborhood basis (in $Y$) $U_1\supseteq U_2\supseteq U_3\supseteq \cdots$ at $y$ and let $0\neq x\in \textnormal{GK}(T)$.
 Then there is some positive integer $p$ such that $T^px=0$. On the other hand, since
$(T,Y)$ is strongly topologically transitive, there is a sequence $(n_k)_{k\ge 1}$ of non-negative integers and a sequence $(y_k)_{k\ge 1}$ with $y_k\in U_k$
such that $x=T^{n_k}y_k$.
Then $0=T^px=T^{p+n_k}y_k$ ($k\ge 1$) which says that, for all $k\ge 1$,
$y_k$ is a vector in the generalized kernel of $T$. But $y_k\to y$ in $Y$ (which implies that $y_k\to y$ in $\mathcal X$), and hence, $y\in\overline{\textnormal{GK}(T)\cap Y}$. 
\end{proof}
\end{prop}
Note that $\overline{\textnormal{GK}(T)\cap Y}\supseteq Y$ if and only if $\overline{\textnormal{GK}(T)\cap Y}=\overline Y$. 
Now we give the following criteria for strong topological transitivity and hypermixing properties which are analogous to those given for bounded linear operators 
\cite[Theorems $2.1$ and $2.3$]{a2}. Since there is no trace of using continuity in the proofs
of those theorems, the similar proofs can be given for the following criteria, and hence, we omit the proofs.
\begin{thm} (STT criterion) Suppose that $T$ is a linear operator in $\mathcal X$
and $S: Y\to Y$ is a map such that $TS=I$ on $Y$. Then $(T,Y)$ is strongly topologically transitive if and only if, for every nonzero $x\in Y$ 
and all $y\in Y$, 
there exist sequences $(n_k)_k$ of positive integers and $(w_k)_k$ in $\textnormal{GK}(T)\cap Y$, with $w_k\in \textnormal{Ker} T^{n_k}$, such that $S^{n_k}x+w_k\to y$ as $k\to \infty$.
\end{thm}
\begin{thm} (HM criterion) Suppose that $T$ is a linear operator in $\mathcal X$
and $S: Y\to Y$ is a map such that $TS=I$ on $Y$. Then $(T,Y)$ is hypermixing if and only if, for every pair of vectors $x, y\in Y$, 
there is a sequence $(w_n)$ in $\textnormal{GK}(T)\cap Y$, with $w_n\in \textnormal{Ker} T^n$, such that $S^nx+w_n\to y$ as $n\to \infty$.
\end{thm}
We will use the following corollary to prove Theorem $3.1$. Again, the proof is omitted since it is similar to that of \cite[Corollary $2.4$]{a2}. 
\begin{cor} Suppose that $T$ is a linear operator in $\mathcal X$
and $S: Y\to Y$ is a map such that $TS=I$ on $Y$. If $\overline{\textnormal{GK}(T)\cap Y}\supseteq Y$ and, for any $x\in Y$,
the sequence $(S^nx)$ is convergent in $Y$ then $(T,Y)$ is hypermixing.
\end{cor}
\subsection{Up-to-closure extensions of a dynamical domain}
In this subsection, we investigate the impact of a certain kind of extension of the dynamical domain of a linear operator on its dynamical properties.
Suppose that $Y_1$ is a dynamical domain for a linear operator $T$ in $\mathcal X$. We say that a dynamical domain 
$Y_2$ for $T$ is an {\it up-to-closure extension} of $Y_1$  
provided that $Y_1\subseteq Y_2\subseteq\overline{Y_1}$. It is clear that if $Y_2$ is an up-to-closure extension of $Y_1$ then $\overline{Y_1}=\overline{Y_2}$.
\par Let P stand for each of the dynamical properties mentioned in this paper except hypercyclicity, strong topological transitivity, and hypermixing. Then we have:
\begin{prop} Suppose $T$ is a linear operator in $\mathcal X$ and $Y_2$ is an up-to-closure extension of $Y_1$.
If $(T,Y_1)$ has the property \textnormal{P} then $(T,Y_2)$ has that property. Also, if $(T,Y_1)$ is hypermixing (resp. strongly topologically transitive) then $(T,Y_2)$ is supermixing
(resp. topologically transitive).
\begin{proof} The key point of the proof is the fact that if $Y_2$ is an up-to-closure extension of $Y_1$ then every nonempty open set in $Y_2$ intersects $Y_1$.
Thus, if $U_2$ is a nonempty open set in $Y_2$ then $U_1=U_2\cap Y_1$ is a nonempty open set in $Y_1$. Now one can easily verify that if $(T,Y_1)$ has the property P then $(T,Y_2)$ has that property. 
\par To prove the second assertion, assume that $(T,Y_1)$ is hypermixing (resp. strongly topologically transitive). Then $(T,Y_1)$ is supermixing (resp. topologically transitive), and hence,
$(T,Y_2)$ is also supermixing (resp. topologically transitive) by the proof of the first part.
\end{proof}
\end{prop}
 Two natural particular cases of up-to-closure extensions of $Y_1=Y$ are $Y_2=\overline Y$ provided that $\overline Y\subseteq Y_0$, and $Y_2=Y_0$ whenever $\overline Y=\overline{Y_0}$. Thus the following results are consequences of Proposition $2.16$.
 \begin{cor} Suppose $T$ is a linear operator in $\mathcal X$ and $\overline Y\subseteq Y_0$. If $(T,Y)$ has the property \textnormal P then $(T,\overline Y)$ has that property.
 Also, if $(T,Y)$ is hypermixing (resp. strongly topologically transitive) then $(T,\overline Y)$ is supermixing (resp. topologically transitive).
 \end{cor}
 If $T$ is bounded and $\mathcal M$ is a closed subspace 
 of $\mathcal X$ then, instead of saying that $(T,\mathcal M)$ is relatively hypercyclic, we can say that $T$ is subspace hypercyclic for $\mathcal M$.
 Thus, as a consequence of Corollary $2.17$, we give the following result.
 \begin{cor} Suppose $T$ is a bounded linear operator on $\mathcal X$. If $(T,Y)$ is relatively hypercyclic then $T$ is subspace hypercyclic for $\overline Y$.
 \end{cor}
 \begin{cor} Suppose $T$ is a linear operator in $\mathcal X$ and $\overline Y=\overline{Y_0}$. If $(T,Y)$ has the property \textnormal P then $T$ has that property.
 Also, if $(T,Y)$ is hypermixing (resp. strongly topologically transitive) then $T$ is supermixing (resp. topologically transitive).
 \end{cor}
 A particular case of the condition $\overline Y=\overline{Y_0}$ in the statement of Corollary $2.19$ is the case that $\overline Y=\mathcal X$.
\begin{cor} Suppose $T$ is a linear operator in $\mathcal X$ and $\overline Y=\mathcal X$. If $(T,Y)$ has the property \textnormal P then $T$ has that property.
 Also, if $(T,Y)$ is hypermixing (resp. strongly topologically transitive) then $T$ is supermixing (resp. topologically transitive).
 \end{cor}
\section{Some hypermixing unbounded linear operators}
In this section, we show that the unbounded linear operators in \cite{bcs} and \cite{jm} which have been proved  to be chaotic, are hypermixing.
In fact, we prove that the derivative operator in $H^2$ and all unbounded weighted left translations in $L_p(0,\infty)$ ($1\le p<\infty$) and $C_0[0,\infty)$ 
are hypermixing (with respect to $Y_0$) and the Laplacian operator in $L^2(\Omega)$ is hypermixing with respect to a proper subspace $Y$ of $Y_0$.
Recall that, if $\Omega$ is a bounded open subset of $\Bbb R^2$, then the Hilbert
space $L^2(\Omega)$ is defined by $$L^2(\Omega)=\{f: \Omega\to \Bbb C: \|f\|^2=\iint_{\Omega}|f(x,y)|^2dxdy<\infty\}.$$
\begin{thm} The following unbounded linear operators are hypermixing:\\
 \hspace*{.2cm} \textnormal{i)} The derivative operator $T=D$ in $H^2$ defined by $Tf=f'$, where
  $Y_0=\bigcap_{n=1}^{\infty}D(T^n)=\{f\in H^2: f^{(n)}\in H^2, \textnormal{for all}\;\;n\ge 1\}$,\\ 
\hspace*{.08cm} \textnormal{ii)} The Laplacian operator $T=(\Delta,Y)$ in $L^2(\Omega)$ defined by 
$Tf=f_{xx}+f_{yy}$, where $Y=\{f\in L^2(\Omega): f \;\;\textnormal{a polynomial}\}$, \\
 \hspace*{.05cm}\textnormal{iii)} All weighted left translations $T=T_{w,a}$ in $L_p(0,\infty)$ ($1\le p<\infty$, $w>1$, $a>0$) 
  defined by $(Tf)(t)=w^tf(t+a)$ ($t>0$), where $Y_0=\bigcap_{n=1}^{\infty}D(T^n)$ and
   $$D(T^n)=\{f\in L_p(0,\infty): \int_{0}^{\infty}|w^{nt+\frac{(n-1)na}{2}}f(t+na)|^pdt<\infty\}\;\;(n\ge 1),$$
   and\\
 \hspace*{.04cm}\textnormal{iv)} All weighted left translations $T=T_{w,a}$ in $C_0[0,\infty)$ ($w>1$, $a>0$)
 defined by $(Tf)(t)=w^tf(t+a)$ ($t\ge 0$), where $Y_0=\bigcap_{n=1}^{\infty}D(T^n)$ and
 $$D(T^n)=\{f\in C_0[0,\infty): \lim_{t\to \infty} w^{nt+\frac{(n-1)na}{2}}f(t+na)=0\}\;\;(n\ge 1).$$ 
\begin{proof} \textnormal{i)}. Let $S$ be the integration operator in $H^2$ defined by $$(Sf)(z)=\int_0^zf(w)dw\;\; (z\in \Bbb D).$$
 Then $TS=I$ on $Y_0$, and it is easy to show that $S$ is in fact a bounded linear operator defined on the entire space $H^2$ satisfying $\|S^n\|\le 1/n$ for all $n\ge 1$. Hence
 $S^nf\to 0$, as $n\to\infty$, for all $f\in Y_0$. On the other hand, it is clear that $\overline{\textnormal{GK}(T)\cap Y_0}=\overline{\textnormal{GK}(T)}=H^2\supseteq Y_0$. Therefore, $T$ is hypermixing
by Corollary $2.15$.\\
\hspace*{.5cm}  \textnormal{ii)} For the right inverse $S=\Delta^{-1}: Y\to Y$ (introduced before \cite[Lemma $3.2$]{bcs}) defined on
$\textnormal{span}\{X_nY_l: n, l\ge 0\}=\textnormal{span}\{x^ny^l: n,l\ge 0\}=Y$
 by
$$\Delta^{-1}(X_nY_l)=\sum_{j=0}^{\infty}(-1)^jX_{n-2j}Y_{l+2j+2},$$
where $X_n=x^n/n!$ and $Y_l=y^l/l!$, we have $S^nf\to 0$, as $n\to\infty$, for every $f\in Y$ \cite[Corollary $3.3$]{bcs}.
On the other hand, $\overline{\textnormal{GK}(T)\cap Y}=\overline{\textnormal{GK}(T)}=L^2(\Omega)\supseteq Y$. Hence $T=(\Delta,Y)$ is hypermixing
by Corollary $2.15$.\\
\hspace*{.5cm}  \textnormal{iii)} Fix an integer $p\ge 1$. For any function $f\in Y_0$, let us define the function $Sf$ in $(0,\infty)$ by 
\begin{equation*}
 (Sf)(t) = \begin{cases}
 w^{-(t-a)}f(t-a)     & t>a,\\
  0             & \text{otherwise}.
\end{cases}
\end{equation*}
We claim that $S(Y_0)\subseteq Y_0$, or equivalently, $S$ maps $Y_0$ into itself. Suppose $f\in Y_0$ is arbitrary. To show that $Sf\in Y_0$, we need to prove that
 $Sf\in D(T^n)$ for all $n\ge 1$, or equivalently, we need to show that $Sf\in L_p(0,\infty)$ and $\|T^n(Sf)\|_p<\infty$ for all $n\ge 1$.
  Indeed, an easy integral calculation shows that $\|Sf\|_p\le \|f\|_p<\infty$, and hence,
$Sf\in L_p(0,\infty)$. Moreover, it is easily seen that $T(Sf)=f$ which says that ($\|T(Sf)\|_p<\infty$ and hence) $Sf\in D(T)$. Now, for all $n\ge 2$, since $f\in D(T^{n-1})$ we have 
$\|T^n(Sf)\|_p=\|T^{n-1}(TSf)\|_p=\|T^{n-1}f\|_p<\infty$, and hence, $Sf\in D(T^n)$. Therefore, the map $S: Y_0\to Y_0$ is a right inverse map for $T$.
On the other hand, we have that $S^nf\to 0$, as $n\to\infty$, for all $f\in Y_0$. In fact, in the proof of \cite[Theorem $4.1$]{jm} (see $4.9$ and $4.10$),
 the map $S$ is defined on $Y=\textnormal{GK}(T)$ and then it is proved that $S^nf\to 0$ for all $f\in Y$ (two lines after $4.12$), but the proof does not
 use the assumption $f\in Y$ (in fact the authors choose $Y=\textnormal{GK}(T)$ because they also need to have $T^nf\to 0$ for all $f\in Y$
  so that they can use Theorem $1.2$). On the other hand, as mentioned in the proof of \cite[Lemma $4.1$]{jm} (see $4.3$ and $4.4$), 
 $\overline {\textnormal{GK}(T)}=L_p(0,\infty)$ and therefore,
  $\overline {\textnormal{GK}(T)\cap Y_0}=\overline {\textnormal{GK}(T)}=L_p(0,\infty)\supseteq Y_0$. Hence $T$ is hypermixing
by Corollary $2.15$.\\
\hspace*{.5cm}  \textnormal{iv)} A similar proof to that given in the previous part shows that, if for any $f\in Y_0$, the function $Sf$ in $[0,\infty)$ is defined by
\begin{equation*}
 (Sf)(t) = \begin{cases}
 \frac{f(0)}{a}t     & 0\le t<a,\\
     w^{-(t-a)}f(t-a)          & t\ge a,
\end{cases}
\end{equation*}
then $S$ maps $Y_0$ into itself and $TSf=f$ for all $f\in Y_0$. Also, the map $S$ satisfies $S^nf\to 0$, as $n\to \infty$, for all $f\in Y_0$. 
Again, in the proof of \cite[Theorem $6.1$]{jm} (see $6.8$ and $6.9$), the map $S$ is defined on $Y=\textnormal{GK}(T)$ and then 
it is proved that $S^nf\to 0$ for all $f\in Y$ (two lines after $6.11$),
   while nowhere in the proof the assumption $f\in Y$ is used. On the other hand,
 as mentioned in the proof of \cite[Lemma $6.1$]{jm} (see $6.3$ and $6.4$), $\overline {\textnormal{GK}(T)}=C_0[0,\infty)$.
Then $\overline{\textnormal{GK}(T)\cap Y_0}=\overline{\textnormal{GK}(T)}=C_0[0,\infty)\supseteq Y_0$. Thus $T$ is hypermixing
by Corollary $2.15$.
\end{proof}
\end{thm}
The following result is a consequence of Theorem $3.1$ and Corollary $2.20$.
\begin{cor} The Laplacian operator $\Delta$ in $L^2(\Omega)$ is supermixing.
\end{cor}
As an application of Theorem $3.1$, we give the following proposition. 
\begin{prop} Let $Y_0=\{f\in H^2: f^{(n)}\in H^2, \textnormal{\;for all}\;\;n\ge 1\}$ and $U$ be any nonempty open subset of $Y_0$. 
Then, for any $0\neq \alpha\in \Bbb C$, there is
a non-negative integer $N$ and a sequence $(p_n)_{n\ge N}$ of polynomials such that, for all $n\ge N$, we have that $p_n\in U$, $\textnormal{deg} p_n=n$, and 
the leading coefficient of $p_n$ is $\alpha/n!$. Furthermore, the assertion is also true for nonempty open sets in the whole space $H^2$.
\begin{proof} Fix a nonempty open set $U$ in $Y_0$ and a complex number $\alpha\neq 0$.
Since $D=(D, Y_0)$ is hypermixing in $H^2$ by Theorem $3.1$ (i), for the constant function $h\in Y_0$ defined by $h(z)= \alpha$ ($z\in \Bbb D$), 
we have that $h\in \bigcap_{i=N}^{\infty}D^n(U)$ for some non-negative integer $N$. Thus, for every integer $n\ge N$, there is a function $p_n\in U$
such that $h=D^np_n$, or equivalently, $\alpha=h(z)=p_n^{(n)}(z)$ for all $z\in \Bbb D$. Then it is obvious that $p_n$ ($n\ge N$) must be a polynomial of degree $n$, and moreover, if $a_n$ is the leading coefficient
of $p_n$, we must have $\alpha=n!a_n$.
\par To prove the second claim, let $G$ be any nonempty open set in $H^2$. Since $\overline Y_0=H^2$ (because $Y_0$ includes the set of all polynomials in $H^2$), 
we have that $G\cap Y_0\neq \emptyset$. Now, if we put $V=G\cap Y_0$ then $V$ is open in $Y_0$, and hence, by the proof of the first part, for any $0\neq \alpha\in \Bbb C$,
there exists a sequence of polynomials in $V$ with the mentioned properties. But $V\subseteq G$ and we are done.
\end{proof}
\end{prop}

\vspace*{.5cm}
\hspace*{.35cm} Department of Mathematics\\
\hspace*{.35cm} Azad University of Gachsaran\\
\hspace*{.35cm} Gachsaran, Iran\\
\hspace*{.35cm} {\it E-mail address}: \email{\small
ansari.moh@gmail.com}, \email{\small mo.ansari@iau.ac.ir}

\end{document}